\documentclass[12pt,reqno]{amsart}

\usepackage{amssymb}
\usepackage{hyperref}

\newcommand{\F}{\mathbb{F}}
\newcommand{\Fp}{\F_p}
\newcommand{\Gal}{\text{\rm Gal}}
\newcommand{\Ic}{\mathcal{I}}
\newcommand{\N}{\mathbb{N}}
\newcommand{\Z}{\mathbb{Z}}

\begin{document}

\title[Galois Module Structure of Milnor $K$-theory]
{Galois Module Structure of Milnor $K$-theory in Characteristic $p$}

\author[Bhandari]{Ganesh Bhandari}
\address{Department of Mathematics, Physics and Engineering, Mount Royal College, \
4825 Mount Royal Gate SW, Calgary, Alberta \ T3E 6K6 \ CANADA}
\email{bhandari@mtroyal.ca}

\author[Lemire]{Nicole Lemire$^\star$}
\address{Department of Mathematics, Middlesex College, \
University of Western Ontario, London, Ontario \ N6A 5B7 \ CANADA}
\thanks{$^\star$Research supported in part by NSERC grant R3276A01.}
\email{nlemire@uwo.ca}

\author[Min\'{a}\v{c}]{J\'an Min\'a\v{c}$^\dag$}
\thanks{$^\dag$Research supported in part by NSERC grant R0370A01, and by a Distinguished Research Professorship at the University of Western Ontario.}
\email{minac@uwo.ca}

\author[Swallow]{John Swallow$^\ddag$}
\address{Department of Mathematics, Davidson College, Box 7046,
Davidson, North Carolina \ 28035-7046 \ USA}
\thanks{$^\ddag$Research supported in part by NSA grant MDA904-02-1-0061 and by NSF grant DMS-0600122.}
\email{joswallow@davidson.edu} 

\begin{abstract}
    Let $E$ be a cyclic extension of degree $p^n$ of a field $F$ of characteristic $p$. Using arithmetic invariants of $E/F$ we determine $k_m E$, the Milnor $K$-groups $K_m E$ modulo $p$, as $\Fp[\Gal(E/F)]$-modules for all $m \in \N$. In particular, we show that each indecomposable summand of $k_m E$ has $\Fp$-dimension a power of $p$.  That all powers $p^i$, $i=0, 1, \dots, n$, occur for suitable examples is shown in a subsequent paper \cite{MSS2}, where additionally the main result of this paper becomes an essential induction step in the determination of $K_m E/p^s K_m E$ as $(\Z/p^s\Z)[\Gal(E/F)]$-modules for all $m, s \in \N$.
\end{abstract}

\date{September 5, 2007}

\maketitle

\newtheorem{theorem}{Theorem}
\newtheorem*{proposition*}{Proposition}
\newtheorem{lemma}{Lemma}
\newtheorem{prop}{Proposition}
\newtheorem*{cor*}{Corollary}
\theoremstyle{definition}
\newtheorem*{remark*}{Remark}
\newtheorem*{reduction}{Reduction to Lemma~\ref{le:normfilt}}

\parskip=10pt plus 2pt minus 2pt

In the 1960s Z.~I.~Borevi\v{c} and D.~K.~Faddeev classified the possible $\Fp[\Gal(E/F)]$-modules $k_1E$, where $E$ is a cyclic extension of degree $p^n$ of a local field $F$. (See \cite{Bo}.) This result has recently been generalized: first in the $n=1$ case, replacing the local field $F$  by any field containing a primitive $p$th-root of unity \cite{MS}; and later  to all base fields $F$ and all cyclic extensions $E$ of degree $p^n$ \cite{MSS1}.  Then three of the authors obtained a generalization from $k_1$ to $k_m$, handling all $\Fp[\Gal(E/F)]$-modules $k_mE$, $m\in \N$, provided that $n=1$ and $F$ contains a primitive $p$th-root of unity \cite{LMS}. This result, in turn, led to new restrictions on possible absolute Galois pro-$p$ groups \cite{BLMS}.

The particular case of fields of characteristic $p$ is important not only for its intrinsic value but also because the case arises naturally in considering the residue fields of valuations on fields of characteristic $0$.  This case, however, has been treated only as far as $k_1 E$ for $E/F$ cyclic of degree $p^n$ \cite{MSS1}.  Twenty years ago Bloch, Gabber, and Kato established an isomorphism between, on one hand, Milnor $K$-theory modulo $p$ in characteristic $p$ and, on the other, the kernel of the Artin-Schreier operator defined on the exterior algebra on K\"ahler differentials \cite{BK}. Not long afterward, two significant papers by Izhboldin considered Milnor $K$-theory in characteristic $p$ and, using the previous work of Bloch, Gabber, and Kato, ascertained its important Galois module properties \cite{I1, I2}.

In this paper we establish that these two results of Izhboldin are enough to determine precisely the Galois module structure of the Milnor groups $K_mE$ mod $p$ when $E$ is a cyclic extension of $p$th-power degree of a field $F$ of characteristic $p$. The result is simple and useful: these modules are direct sums of trivial modules and modules free over some quotient of the Galois group. Equivalently, the dimensions over $\Fp$ of the indecomposable $\Fp[\Gal(E/F)]$-modules occurring as direct summands of $K_m E/pK_m E$ are all powers of $p$. In \cite{MSS2} it is shown that in fact each indecomposable summand of dimension $p^i$, $i=0,1,\dots,n$, occurs in suitable examples.

Our decomposition of the $\Fp[\Gal(E/F)]$-module $k_mE$ into indecomposables in Theorem~\ref{th:main} is not canonical, although the ranks of the summands which are free over the various quotients of $\Gal(E/F)$ are invariants of the $\Fp[\Gal(E/F)]$-module $k_mE$. Still, our results do have a canonical formulation, stemming from Lemma~\ref{le:normfilt}, which is in fact equivalent to Theorem~\ref{th:main}.  In particular, the filtration of $(k_mE)^{\Gal(E/F)}\simeq k_mF$, induced by the images of successive powers of the maximal ideal of $\Fp[\Gal(E/F)]$, is determined by the images of the norm maps $k_mE_i\to k_mF$, where $E_i$ runs through the intermediate fields between $E$ and $F$.

It is possible that one could use the results in \cite{J} to obtain alternative proofs, but we preferred a more self-contained, short derivation instead.

In \cite{MSS2} our main result is used as a critical first induction step to determine $K_m E/p^s K_m E$ as a $(\Z/p^s\Z) [\Gal(E/F)]$-module for all $m, s \in \N$. The classification problem of $(\Z/p^s\Z) G$-modules for cyclic $G$ is nontrivial and in fact still not complete, despite some important results when, for instance, the module is free over $\Z/p^s\Z$. (See \cite{T} and the references contained therein.)  However, one can describe the $(\Z/p^s\Z) [\Gal(E/F)]$-module $K_m E/p^s K_m E$ completely and simply.

We assume in what follows that all fields have characteristic $p$ and that $m$ is a fixed natural number. For a field $F$, let $K_mF$ denote the $m$th Milnor $K$-group of $F$ and $k_mF=K_mF/ pK_mF$. (See, for instance, \cite{Ma} and \cite[IX.1]{Mi}.) If $E/F$ is a Galois extension of fields, let $G=\Gal(E/F)$ denote the associated Galois group. When $G$ is a cyclic group we write $G = \langle \sigma \rangle$, with a suitable fixed generator $\sigma$. For the sake of simplicity we write $\rho$ instead of $\sigma-1$. We write $i_E\colon K_mF\to K_mE$ and $N_{E/F}\colon K_mE\to K_mF$ for the natural map induced by the inclusion and the norm maps, and we use the same notation for the induced maps modulo $p$ between $k_mF$ and $k_mE$. In order to avoid possible confusion, in a few instances we write $i_{F,E}$ instead of $i_E$.

Izhboldin's results are as follows.

\begin{lemma}[{\cite[Lemma~2.3]{I2}}]\label{le:i1}
    Suppose $E/F$ is cyclic of degree $p$.  Then $i_E\colon
    k_mF\to (k_mE)^G$ is an isomorphism.
\end{lemma}

\begin{theorem}[Hilbert 90 for Milnor $K$-theory: {\cite[Corollary of
Proposition~5]{I1}}, {\cite[Theorem~D]{I2}}]\label{th:exseq}
    Suppose $E/F$ is cyclic of $p$th-power degree. Then the
    following sequence is exact:
    \begin{equation*}
        K_mE \xrightarrow{1-\sigma} K_mE \xrightarrow{N_{E/F}}
        K_mF.
    \end{equation*}
\end{theorem}

Now suppose $E/F$ is cyclic of degree $p^n$.  For $i=0, 1, \dots, n$, let $E_i/F$ be the subextension of degree $p^i$ of $E/F$, $H_i=\Gal(E/E_i)$, and $G_i=\Gal(E_i/F)$. Our main result is the following
\begin{theorem}\label{th:main}
    There exists an isomorphism of $\Fp G$-modules
    $k_mE \simeq \oplus_{i=0}^n Y_i$, where
    \begin{itemize}
    \item $Y_n$ is a free $\Fp G$-module of rank $\dim_{\Fp}
    N_{E/F}k_mE$,
    \item $Y_i$, $0<i<n$, is a free $\Fp G_i$-module of rank
    \begin{equation*}
        \dim_{\Fp} N_{E_i/F} k_mE_i/N_{E_{i+1}/F}k_mE_{i+1},
    \mbox{ and }
    \end{equation*}
    \item $Y_0$ is a trivial $\Fp G$-module of rank $\dim_{\Fp}
    k_mF/N_{E_1/F} k_mE_1$.
    \end{itemize}
\end{theorem}

Our proof of Theorem~\ref{th:main} relies on a filtration lemma, Lemma~\ref{le:normfilt}, which shows that certain elements in $k_mE$ are expressible as norms.  That the various norm groups $N_{E_i/F}k_mE_i$ contain enough elements is, in fact, the key to this lemma.  In section~\ref{se:n1} we present technical results on $\Fp G$-modules. In section~\ref{se:n2} we state Lemma~\ref{le:normfilt}, prove that Theorem~\ref{th:main} is a corollary to Lemma~\ref{le:normfilt}, and prove basic lemmas in preparation for its proof.  This proof is finally presented in section~\ref{se:n3}.

\section{$\Fp G$-modules}\label{se:n1}

For the reader's convenience, after introducing some notation,
we recall in this section some basic elementary facts about
$\Fp G$-modules. (See \cite{Ca} and \cite{L}.)

Let $G$ be a cyclic group of order $p^n$ with generator $\sigma$. We denote $\sigma-1$ as $\rho$. For an $\Fp G$-module $U$, let $U^G$ denote the submodule of $U$ fixed by $G$.  For an arbitrary element $u\in U$, we say that the length $l(u)$ is the dimension over $\Fp$ of the $\Fp G$-submodule $\langle u\rangle$ of $U$ generated by $u$. Then we have
\begin{equation*}
    \rho^{l(u)-1}\langle u \rangle =
    \langle u \rangle^G \neq \{0\} \text{\quad and\quad}
    \rho^{l(u)}\langle u \rangle = \{0\}.
\end{equation*}
If we want to stress the dependence of length $l(u)$ on $G$, we write $l_G(u)$. As usual, if $U$ is a free $\Fp G$-module with $U = \oplus_{i\in\Ic} \Fp G$, we say that $U$ is a module of rank $\vert\Ic\vert$. Denote by $N$ the operator $\rho^{p^n-1}$ acting on module $U$.

We use the following general lemma about $\Fp G$-modules. The simple proof of this lemma is omitted.

\begin{lemma}[Exclusion Lemma]\label{le:excl}
    Let $M_i$, $i \in \Ic$, be a family of $\Fp G$-modules
    contained in a common $\Fp G$-module $N$. Suppose that
    the $\Fp$-vector subspace $R$ of $N$ generated by all
    $M_i^G$ has the form $R = \oplus_{i \in \Ic} M_i^G$. Then
    the $\Fp G$-module $M$ generated by $M_i$, $i \in \Ic$,
    has the form $M = \oplus_{i \in \Ic} M_i$.
\end{lemma}

We also use the following result about the restriction of a
cyclic $\Fp G$-module to a subgroup.

\begin{lemma}[Restriction Lemma]\label{le:rest}
    Let $Y$ be a cyclic $\Fp G$-module of length $l$ generated
    by $\gamma$. Then $l_{H_k}(\gamma)$ is the unique integer
    such that $l=(l_{H_k}(\gamma)-1)p^k+r$ for some unique integer
    $1\le r\le p^k$. Moreover,
    \begin{equation*}
        Y\simeq V_{l_{H_k}(\gamma)}^{r}\oplus V_{l_{H_k}(\gamma)-1}^{p^k-r}
    \end{equation*}
    where $V_i$ is a cyclic $\Fp H_k$-module of length $i$.
\end{lemma}

\begin{proof}
    The formula for $l_{H_k}$ follows immediately from the definition of length. Since $\{\rho^i\gamma\ \vert\  0\le i\le l-1\}$ is an $\Fp$-basis for $Y=\Fp G\gamma$, then
    \begin{equation*}
        \Fp G\gamma\big\vert_{H_k}\ =\ \oplus_{i=0}^{p^k-1}\Fp H_k
        \rho^i\gamma \ \simeq\   V_{l_{H_k}(\gamma)}^{r}\oplus V_{l_{H_k}
        (\gamma)-1}^{p^k-r}
    \end{equation*}
    since $l_{H_k}(\rho^j\gamma)=l_{H_k}(\gamma)$ if $0 \le j < r$ and $l_{H_k} (\rho^j\gamma)=l_{H_k}(\gamma)-1$ if $r \le j \le p^k-1$.
\end{proof}

Finally, we need a general structure proposition about $\Fp G$-modules that shows that the structure of an $\Fp G$-module $X$ can be determined from the natural filtration on $X^G$ obtained by taking the intersection of $X^G$ with the images of $X$ under the successive powers of the augmentation ideal of $\Fp G$. The proposition below is proved for cyclic $p$-groups as \cite[Proposition~2]{LMS}.  The generalization to cyclic groups of order $p^n$ is automatic.

\begin{prop}\label{pr:fpgstruct}
    Let $X$ be an $\Fp G$-module. Set $L_{p^n}=\rho^{p^n-1}X,$ and for $1\le i<p^n$, suppose that $L_i$ is an $\Fp$-complement of $\rho^{i}X\cap X^G$ in $\rho^{i-1}X\cap X^G$.

    Then there exist $\Fp G$-modules $X_i$, $i=1, 2, \dots, p^n$, such that
    \begin{enumerate}
        \item $X=\bigoplus_{i=1}^{p^n} X_i$,
        \item $X_i^G=L_i$ for $i=1, 2, \dots, p^n$,
        \item each $X_i$ is a direct sum of $\dim_{\Fp}(L_i)$ cyclic
        $\Fp G$-modules of length $i$, and
        \item for each $i=1, 2, \dots, p^n$, there exists an $\Fp$-submodule
        $Y_i$ of $X_i$ with $\dim_{\Fp}(Y_i)= \dim_{\Fp}(L_i)$ such that
        $\Fp GY_i=X_i$.
    \end{enumerate}
\end{prop}

\section{Milnor $k$-groups}\label{se:n2}

If $\alpha\in K_mE$ we write $\bar\alpha$ for the class of $\alpha$ in $k_mE$.  For $\gamma\in K_mE$, let $l(\gamma)$ denote the dimension over $\Fp$ of the $\Fp G$-submodule $\langle \bar\gamma \rangle$ of $k_mE$ generated by $\bar\gamma$. Because $\rho^{p^n-1} = 1+\sigma+\dots+\sigma^{p^n-1}$ in $\Fp G$, we may use $i_EN_{E/F}$ and $N$ interchangeably on $k_mE$.

We will prove the following key lemma by induction on $n$.

\begin{lemma}~\label{le:normfilt}
    Let $1\le j\le p^n$ and let $0\le i\le n$ be minimal such that $j\le p^i$. Then
    \begin{equation*}
        \rho^{j-1}k_mE\cap (k_mE)^G=i_EN_{E_i/F}k_mE_i.
    \end{equation*}
\end{lemma}

\begin{remark*}
    Note that Lemma~\ref{le:i1} proves this result for $j=1$ and $n=1$. Note also that for $j=p^n$, the result holds by definition. Moreover, for $i$ minimal such that $j\le p^i$,
    \begin{align*}
        i_EN_{E_i/F}k_mE_i&= i_{E_i,E}i_{F,E_i}N_{E_i/F}k_mE_i =
        \rho^{p^i-1}(i_{E_i,E}k_mE_i)\\ &\subseteq \rho^{p^i-1}
        k_mE\cap (k_mE)^G \subseteq \rho^{j-1}k_mE\cap (k_mE)^G,
    \end{align*}
    so it suffices to show the inclusion in the opposite
    direction.
\end{remark*}

Before proving this lemma, we show
\begin{reduction}\label{le:reduc}
    Theorem~\ref{th:main} holds for those cyclic extensions of degree $p^n$ for which Lemma~\ref{le:normfilt} holds.
\end{reduction}

\begin{proof}
    We apply Proposition~\ref{pr:fpgstruct} to $k_mE$. Set $L_{p^n} = i_EN_{E/F}k_mE$, and for $0\le i<n$, define $L_{p^i}$ as an $\Fp$-complement of $i_E N_{E_{i+1}/F}k_mE_{i+1}$ in $i_EN_{E_i/F}k_mE_i$.  We also set $L_k = \{0\}$ for each $k \in \{ 1, 2, \dots, p^n \}$ which is not a power of $p$.  Note that $(k_mE)^G= i_Ek_mF$ by Lemma~\ref{le:normfilt} applied to $j=1$. Then, Lemma~\ref{le:normfilt} shows that the $\Fp$-modules $L_j$, $j=1,\dots,p^n$, satisfy the conditions of the proposition. By the proposition we have a decomposition of $k_mE$ into a direct sum of $\Fp G$-modules
    \begin{equation*}
        k_mE=\oplus_{i=0}^n Y_i,
    \end{equation*}
    where $Y_i^G=L_{p^i}$ and $Y_i$ is a direct sum of
    \begin{align*}
        \dim_{\Fp}L_{p^i} &= \dim_{\Fp}i_EN_{E_{i}/F}k_mE_i/ i_EN_{E_{i+1}/F}k_mE_{i+1}\\ &=\dim_{\Fp}N_{E_i/F}k_mE_i/ N_{E_{i+1}/F}k_mE_{i+1}
    \end{align*}
    cyclic $\Fp G$-modules of length $p^i$, using Lemma~\ref{le:i1} for the last equality.  Equivalently, $Y_i\simeq \oplus_{T_i} \Fp(G/H_i)$, where the cardinal number of $T_i$ is $\dim_{\Fp}L_{p^i}$, and therefore $Y_i$ is a free $\Fp G_i$-module of the same rank.
\end{proof}

We now prove  several lemmas as partial results toward the proof of Lemma~\ref{le:normfilt}. Note that it is only in these lemmas that Theorem~\ref{th:exseq} is used in this paper.

\begin{lemma}\label{le:key}
    Lemma~\ref{le:normfilt} holds for all cyclic extensions  $E/F$ of degree $p>2$.
\end{lemma}

\begin{proof}
    Let $G=\Gal(E/F)$. By the remark to Lemma~\ref{le:normfilt}, it suffices to show the left-to-right inclusion of Lemma~\ref{le:normfilt}. Let $\gamma\in K_mE$ such that $0\ne \rho^{l-1}\bar \gamma\in \rho^{l-1}k_mE\cap (k_mE)^G$ where $1<l\le p$. Then $l(\gamma)=l$.
    We show by induction on $i$, $l\le i\le p$, that there exists $\alpha_i\in K_mE$ such that
    \begin{equation*}
        \langle \rho^{i-1}\bar\alpha_i\rangle = \langle
        \bar\gamma \rangle^G.
    \end{equation*}
    Then since $\rho^{p-1}\bar\alpha_p = i_E N_{E/F} \bar\alpha_p$,the proof will be complete. If $i=l$ then $\alpha_i=\gamma$ suffices.  Assume now that $l\leq i<p$ and our statement is true for $i$.

    Set $c := N_{E/F}\alpha_i$.  Since $i_E\bar c = N\bar\alpha_i
    = \rho^{p-1}\bar\alpha_i$ and $i<p$, we conclude that
    $i_E\bar c = 0$.  By the injectivity of $i_E$ from
    Lemma~\ref{le:i1}, $\bar c = 0$.  Therefore there exists
    $f\in K_mF$ such that $c=pf$ in $K_mF$.

    We calculate
    \begin{equation*}
        N_{E/F}(\alpha_i-i_E(f)) = c-pf = 0.
    \end{equation*}
    By Theorem~\ref{th:exseq}, there exists $\omega\in K_mE$ such
    that $\rho\omega = \alpha_i - i_E(f)$. Hence $\rho^2\omega =
    \rho\alpha_i$.  Since $i\ge 2$,
    \begin{equation*}
        \langle \rho^i \bar\omega \rangle = \langle
        \rho^{i-1}\bar\alpha_i\rangle = \langle \bar\gamma
        \rangle^G
    \end{equation*}
    and we can set $\alpha_{i+1}=\omega$.
\end{proof}

\begin{lemma}\label{le:key2}
    Lemma~\ref{le:normfilt} holds for all cyclic extensions $E/F$ of degree $2$ or $4$.
\end{lemma}

\begin{proof}
    By the remark to Lemma~\ref{le:normfilt}, it suffices to prove that for $E/F$ a cyclic extension of degree $4$,
    \begin{equation*}
        \rho^2 k_mE\cap (k_mE)^G\subseteq i_EN_{E/F}k_mE.
    \end{equation*}
    Let $\gamma\in K_mE$ such that  $0\ne \rho^2 \bar\gamma\in \rho^2 k_mE\cap (k_mE)^G$. Then $l(\gamma)=3$, and set $\beta := \rho\gamma$. We have
    \begin{equation*}
        (\sigma^2-1)\bar\beta = \rho^3\bar\gamma = 0.
    \end{equation*}
    By Lemma~\ref{le:i1} applied to $E/E_1$, $\bar\beta \in i_{E_1,E} k_mE_1$. Let $b\in K_mE_1$ such that $i_{E_1,E}\bar b = \bar\beta$.

    Set $c := N_{E_1/F} b$. Then
    \begin{equation*}
        i_{F,E}\bar c = i_{E_1,E}(\sigma+1)\bar b =
        (\sigma^2-1)\bar\gamma = i_{E_1,E}N_{E/E_1}\bar\gamma,
    \end{equation*}
    and since $l(\gamma)=3$, $\langle \bar\gamma\rangle^G = \langle
    i_{F,E}\bar c\rangle$.  By Lemma~\ref{le:i1} we see that
    $i_{F,E_1}\bar c= N_{E/E_1}\bar\gamma$ and therefore
    $N_{E/E_1}\gamma = i_{F,E_1}c+2g$ for some $g\in K_mE_1$.

    Now $N_{E/F}\gamma = N_{E_1/F}(i_{F,E_1}c+2g) = 2c + 2N_{E_1/F}g$. Set $\delta := b + g$. We calculate
    \begin{equation*}
        N_{E/F}i_{E_1,E}\delta = 2c+2N_{E_1/F}g,
    \end{equation*}
    so that $N_{E/F}(\gamma-i_{E_1,E}\delta) = 0$. By Theorem~\ref{th:exseq} there exists an $\alpha\in K_mE$ with $\rho\alpha = \gamma-i_{E_1,E}\delta$. Then, observing that $\rho^2 (\overline{i_{E_1,E} \delta}) = 0$ since $\delta \in K_m E_1$, we see that
    \begin{equation*}
        i_{F,E}N_{E/F}\bar\alpha = \rho^3\bar\alpha =
        \rho^2\overline{(\gamma-i_{E_1,E}\delta)} = \rho^2 \bar\gamma,
    \end{equation*}
as required.
\end{proof}

\section{Proof of Lemma~\ref{le:normfilt}}\label{se:n3}

We prove Lemma~\ref{le:normfilt} and hence, using the Reduction to Lemma~\ref{le:normfilt}, Theorem~\ref{th:main}.  We do so by induction on $n$, using Lemmas~\ref{le:key} and \ref{le:key2} as base cases and assuming the result for cyclic extensions of degree $p^{n-1}$. Observe that Lemmas~\ref{le:key} and \ref{le:key2} and the Reduction to Lemma~\ref{le:normfilt} give us that Theorem~\ref{th:main} holds for degree $p$ and, if $p=2$, then for degree $4$ as well. By the remark to Lemma~\ref{le:normfilt}, it suffices to prove that
\begin{equation*}
    \rho^{j-1}k_mE\cap (k_mE)^G \subseteq i_EN_{E/F}k_mE
\end{equation*}
for all $p^{n-1}<j< p^n$.

Assume that $k_mE=\oplus_{\gamma\in \Gamma}\Fp G\gamma$ is a decomposition of $k_mE$ into cyclic $\Fp G$-modules given by Proposition~\ref{pr:fpgstruct}.  Cyclic $\Fp G$-modules are indecomposable, and they are moreover self-dual and local, therefore with local endomorphism rings.  By the Krull-Schmidt-Azumaya Theorem (see \cite[Theorem~12.6]{AF}), all decompositions of $J$ into indecomposables are equivalent. (In our special case one can check this fact directly.)

First assume that $p>2$. We know by Lemma~\ref{le:key} and by the Reduction to Lemma~\ref{le:normfilt} that $k_mE\vert_{H_{n-1}}$ is a direct sum of trivial and free $\Fp H_{n-1}$-modules.  Suppose that $\gamma\in k_mE$ satisfies $p^{n-1}<l(\gamma)< p^n$.  We want to show that $\rho^{l(\gamma)-1}$ $\gamma\in\rho^{p^n-1}k_m E = i_E k_m E$. Without loss of generality we may assume that $\gamma$ generates a cyclic $\Fp G$-summand of $k_mE$.  Since $l(\gamma)<p^n$, then by Lemma~\ref{le:rest}, $\Fp G \gamma\vert_{H_{n-1}}$ contributes direct summands which are not free or trivial, a contradiction. Hence a cyclic generator $\gamma$ with length greater than $p^{n-1}$ must have length $p^n$.

Now assume $p=2$.  Suppose that $\gamma\in \Gamma$ satisfies $2^{n-1}<l(\gamma)<2^{n}$. By Lemma~\ref{le:key2} and the Reduction to Lemma~\ref{le:normfilt}, $k_mE \vert_{H_{n-2}}$ is a direct sum of cyclic $\F_2 H_{n-2}$-modules of lengths 1, 2, and 4, and $\F_2 G \gamma$ is a cyclic summand of $k_mE$.  By Lemma~\ref{le:rest}, we know that $\F_2 G\gamma \vert_{H_{n-2}}$ contributes cyclic $\F_2 H_{n-2} $-summands of length $3$, a contradiction. Hence a cyclic generator $\gamma\in \Gamma$ of length greater than $2^{n-1}$ must have length $2^n$.

In both cases, we have shown for $p^{n-1}<j \le p^n$ that
\begin{equation*}
    \rho^{j-1}k_mE\cap (k_mE)^G = i_E N_{E/F} k_mE,
\end{equation*}
as desired.
\qed

\section*{Acknowledgments}
We are very grateful to A.~Schultz as some of his ideas, developed in \cite{MSS1}, proved quite useful to us during the investigations leading to this paper. We also wish to thank M.~Rost and J.-P.~Serre for their comments after a lecture delivered in Nottingham by the third author in September 2005. These comments led to the improvement of our exposition.

\end{document}